\documentclass[11pt,leqno,a4paper]{amsart}
\usepackage{amsmath,amssymb,amscd,latexsym}

\newcommand{\defeq}{\ensuremath{\,\raisebox{.08ex}{:}\hspace{-.7ex}=}}
\numberwithin{equation}{section}
\usepackage[all]{xy}
\newcommand{\C}{{\mathbb{C}}}
\newcommand{\R}{{\mathbb{R}}}
\newcommand{\Z}{{\mathbb{Z}}}
\newcommand{\ddet}{\mathrm{det}}
\newcommand{\tr}{\mathrm{tr}}
\newcommand{\Nh}{{\mathcal N}}
\newcommand{\halb}{\frac{1}{2}}
\newtheorem{theo}{Theorem}[section] 
\newtheorem{prop}[theo]{Proposition} 
\newtheorem{lemm}[theo]{Lemma} 
\newtheorem{coro}[theo]{Corollary} 

\theoremstyle{definition}
\newtheorem{defi}[theo]{Definition} 
\newtheorem{exam}[theo]{Example} 

\theoremstyle{remark}
\newtheorem{rema}[theo]{Remark}

\begin{document}
\setlength{\baselineskip}{1.1\baselineskip}

\title[Expansive algebraic actions]{Expansive algebraic actions of discrete residually finite amenable groups\\and their entropy}

\author{Christopher Deninger} \address{Christopher Deninger: Mathematical Institute, Westf. Wilhelms-University M\"unster, Einsteinstr. 62, 48149 M\"unster, Germany} \email{c.deninger@math.uni-muenster.de}

\author{Klaus Schmidt} \address{Klaus Schmidt: Mathematics Institute, University of Vienna, Nordberg\-stra{\ss}e 15, A-1090 Vienna, Austria \newline\indent \textup{and} \newline\indent Erwin Schr\"odinger Institute for Mathematical Physics, Boltzmanngasse~9, A-1090 Vienna, Austria} \email{klaus.schmidt@univie.ac.at}

\date{}

	\begin{abstract}
We prove an entropy formula for certain expansive actions of a countable discrete residually finite group $\Gamma $ by automorphisms of compact abelian groups in terms of Fuglede-Kadison determinants. This extends an earlier result proved by the first author under somewhat more restrictive conditions.

The main tools for this generalization are a representation of the $\Gamma $-action by means of a `fundamental homoclinic point' and the description of entropy in terms of the renormalized logarithmic growth-rate of the set of $\Gamma _n$-fixed points, where $(\Gamma _n,\,n\ge1)$ is a decreasing sequence of finite index normal subgroups of $\Gamma $ with trivial intersection.
	\end{abstract}


\thanks{The second author (KS) was partially supported by the FWF Research Grant P16004-N05}

\maketitle

\section{Introduction\label{s:introduction}}

For a countable discrete amenable group $\Gamma$ and an element $f = \sum f_{\gamma} \gamma$ in the integral group ring $\Z \Gamma$ consider the quotient $\Z \Gamma / \Z \Gamma f$ of $\Z\Gamma$ by the left ideal $\Z\Gamma f$ generated by $f$. It is a discrete abelian group with a left $\Gamma$-action by multiplication. The Pontryagin dual
	\begin{equation}
	\label{eq:1.1}
X_f = \widehat{\Z \Gamma / \Z \Gamma f}
	\end{equation}
is a compact abelian group with a left action of $\Gamma$ by continuous group automorphisms. More explicitly, $X_f$ is the following closed subshift of the full shift $(\R / \Z)^{\Gamma} = \widehat{\Z \Gamma}$ with $\Gamma$-action given by $(\gamma \cdot x)_{\gamma'} = x_{\gamma^{-1} \gamma'}$. The elements of $X_f$ are the sequences $(x_{\gamma'})$ in $(\R / \Z)^{\Gamma}$ which satisfy the equations \[ \sum_{\gamma'} f_{\gamma'} x_{\gamma \gamma'} = 0 \quad \mbox{in} \; \R / \Z \; \mbox{for all} \; \gamma \; \mbox{in} \; \Gamma . \] For $\Gamma = \Z^d$ the theory of such dynamical systems has been extensively studied, see \cite{Sch} for example. For general $\Gamma $ the investigation of such systems was initiated in \cite{ER}.

The various definitions of entropy, topological, metric and measure theoretic with respect to Haar measure, all coincide for this action, see e.g. \cite[section 2]{De}. We denote by $h_f$ this common entropy.

If $f$ is a unit in the convolution algebra $L^1 (\Gamma)$, then the action of $\Gamma$ on $X_f$ is expansive. If in addition $f$ is positive in the von~Neumann algebra $\Nh\Gamma$ and $\Gamma$ has a ``log-strong F{\o}lner sequence'', it was shown in \cite{De} that \[ h_f = \log \ddet_{\Nh\Gamma} f , \] where $\ddet_{\Nh\Gamma}$ is the Fuglede--Kadison determinant, \cite{FK}, \cite{L}. Currently the only groups known to have a $\log$-strong F{\o}lner sequence are the virtually nilpotent ones, so this condition is unfortunate, even though it is only needed for the inequality $h_f \le \log \ddet_{\Nh\Gamma} f$ in \cite{De}.

Recall that a countable group $\Gamma$ is residually finite if there is a sequence $\Gamma_n$ of normal subgroups of finite index whose intersection is trivial.

Our main result is this:
	\begin{theo}
	\label{theorem1}
Let $\Gamma$ be a countable discrete amenable and residually finite group and $f$ an element of $\Z \Gamma$. Then the action of $\Gamma$ on $X_f$ is expansive if and only if $f$ is a unit in $L^1 (\Gamma)$. In this case we have the formula
	\begin{equation}
	\label{eq:1.3}
h_f = \log \ddet_{\Nh\Gamma} f .
	\end{equation}
	\end{theo}
For example, $\Gamma$ could be any finitely generated solvable subgroup of a matrix group over a field.

If $f = \sum f_{\gamma} \gamma$ is invertible in $L^1 (\Gamma)$ then $\ddet_{\Nh\Gamma} f =\ddet_{\Nh\Gamma} f ^*$, where $f^* = \sum f_{\gamma} \gamma^{-1}$, since the Fuglede-Kadison determinant of an invertible operator is equal to that of its adjoint. From \eqref{eq:1.3} we therefore obtain that $h_f=h_{f^*}$. This question had been left open in \cite{De} even for the integral Heisenberg group.

Dynamically the proof of formula \eqref{eq:1.3} is based on a description of the entropy for expansive $\Gamma$-actions as a renormalized logarithmic growth rate of the number of $\Gamma_n$-fixed points, where $(\Gamma _n,\,n\ge1)$ is a decreasing sequence of finite index normal subgroups with trivial intersection. In the case of expansive $\Z^d$-actions this relation between entropy and periodic points is known and the methods can be carried over to our case. It remains to prove that $\ddet_{\Nh\Gamma} f$ is the limit of certain finite dimensional determinants. This is a well-known problem in the context of $L^2$-invariants for residually finite groups c.f. \cite[Chapter 13]{L}. There one tries to approximate $L^2$-Betti numbers, $L^2$-signatures etc.~of coverings with covering group $\Gamma$ by the corresponding quantities for the finite groups $\Gamma / \Gamma_n$. For the combinatorial $L^2$-torsion which involves Fuglede--Kadison determinants there is no general result in this direction however.

Nonetheless in our situation which is quite favourable analytically it is not difficult to prove the required approximation property of $\ddet_{\Nh\Gamma} f$ directly.

It is a pleasure for us to thank our colleagues Gabor Elek, Doug Lind and Wolfgang L\"uck for helpful comments. We are also grateful to Doug Lind for alerting us to the reference \cite{Weiss}. The first author would like to thank the Erwin Schr\"odinger Institute for support. 

\section{Group rings
	\label{s:invertibility}
}

Let $\Gamma $ be a countable discrete group with identity element $1=1_{\Gamma}$ and integral group ring $\mathbb{Z}\Gamma$. We denote by $L^\infty (\Gamma )$ the set of all bounded maps from $\Gamma $ to $\mathbb{R}$, write every $w\in L^\infty (\Gamma )$ as $(w_\gamma )$ with $w_\gamma \in\mathbb{R}$ for every $\gamma \in\Gamma $, and denote by $\|w\|_\infty =\sup_{\gamma \in\Gamma }|w_\gamma |$ the supremum norm on $L^\infty (\Gamma )$. For $p\in[1,\infty )$ we set
	$$
L^p(\Gamma )=\biggl\{w\in L^\infty (\Gamma ):\|w\|_p\defeq\biggl(\sum_{\gamma \in\Gamma }|w_\gamma |^p\biggr)^{1/p}<\infty \biggr\}.
	$$

For $\gamma \in\Gamma $ we define $e(\gamma )\in L^1(\Gamma )\subset L^\infty (\Gamma )$ by
	\begin{equation}
	\label{eq:e}
e(\gamma )_{\gamma '}=
	\begin{cases}
1&\textup{if}\enspace \gamma =\gamma '
	\\
0&\textup{otherwise}.
	\end{cases}
	\end{equation}
Every $h$ in $L^1 (\Gamma)$ can be uniquely written as a convergent series
	\begin{equation}
	\label{eq:h}
h=\sum_{\gamma \in\Gamma }h_\gamma e(\gamma )
	\end{equation}
with $h_\gamma \in\mathbb{R}$ and $\sum_{\gamma \in\Gamma }|h_\gamma |<\infty $. In this notation multiplication in the convolution algebra $L^1 (\Gamma)$ takes the form
	\begin{equation}
	\label{eq:product}
	\begin{aligned}
h\cdot h'=\sum_{\gamma ,\gamma '\in\Gamma }h_\gamma h'_{\gamma '}e(\gamma \gamma ')&=\sum_{\gamma \in\Gamma }\biggl(\sum_{\gamma '\in\Gamma }h_{\gamma \gamma '^{-1}}h'_{\gamma '}\biggr)e(\gamma )
	\\
&= \sum_{\gamma \in\Gamma }\biggl(\sum_{\gamma '\in\Gamma }h_{\gamma '}h'_{{\gamma '}^{-1}\gamma }\biggr)e(\gamma ) .
	\end{aligned}
	\end{equation}
The involution $h\mapsto h^*$ in $L^1 (\Gamma)$ is defined by
	\begin{equation}
	\label{eq:involution}
h^*=\sum_{\gamma \in\Gamma }h_{\gamma^{-1}}e(\gamma ) = \sum_{\gamma \in\Gamma }h_\gamma e(\gamma ^{-1})
	\end{equation}
and satisfies that $(g\cdot h)^*=h^*\cdot g^*$ for $g,h\in L^1 (\Gamma)$. An element $h\in L^1(\Gamma )$ is \emph{self-adjoint} if $h^*=h$.

Note that the integral group ring $\Z\Gamma$ of finite formal sums $\sum_{\gamma} a_{\gamma} \gamma$ with $a_{\gamma} \in \Z$ can be viewed as a subring of $L^1 (\Gamma)$ by identifying $\sum_{\gamma} a_{\gamma} \gamma$ with $\sum_{\gamma} a_{\gamma} e (\gamma)$.

Multiplication in $L^1(\Gamma )$ can be extended further to commuting (left) actions $(h,w)\mapsto \lambda_hw$ and $(h,w)\mapsto \rho_hw$ of $L^1(\Gamma )$ on $L^\infty (\Gamma )$ by setting
	\begin{equation}
	\label{eq:action}
	\begin{aligned}
(\lambda_hw)_\gamma &=(h\cdot w)_\gamma =\sum_{\gamma '\in\Gamma }h_{\gamma {\gamma '}^{-1}}w_{\gamma '}=\sum_{\gamma '\in\Gamma }h_{\gamma '}w_{{\gamma '}^{-1}\gamma },
	\\
(\rho_hw)_\gamma &=(w\cdot h^*)_\gamma =\sum_{\gamma '\in\Gamma }w_{\gamma {\gamma '}^{-1}}h^*_{\gamma '}=\sum_{\gamma '\in\Gamma }w_{\gamma \gamma '}h_{\gamma '}
	\end{aligned}
	\end{equation}
for every $h\in L^1(\Gamma )$, $w\in L^\infty (\Gamma )$ and $\gamma \in\Gamma $. The \emph{left} and \emph{right shift actions} $\lambda$ and $\rho$ of $\Gamma $ on $L^\infty (\Gamma )$ are given by
	$$
\lambda^\gamma w=\lambda_{e(\gamma )}w, \qquad \rho^\gamma w=\rho_{e(\gamma)} w
	$$
or, equivalently, by
	\begin{equation}
	\label{eq:shift}
(\lambda^\gamma w)_{\gamma '}=w_{\gamma ^{-1}\gamma '}, \qquad (\rho^\gamma w)_{\gamma '}=w_{\gamma '\gamma }
	\end{equation}
for every $w\in L^\infty (\Gamma )$ and $\gamma ,\gamma '\in\Gamma $.

For the following lemmas we fix $h\in L^1(\Gamma )$ and define the linear operators $\rho_h,\rho_{h^*}\colon L^\infty (\Gamma )\longrightarrow L^\infty (\Gamma )$ by \eqref{eq:involution}--\eqref{eq:action}. For $v\in L^p(\Gamma )$ and $w\in L^q(\Gamma )$ with $\frac1p+\frac1q=1,\,1\le p,q\le\infty $, we set
	\begin{equation}
	\label{eq:innerproduct}
\langle v,w\rangle =\sum_{\gamma \in\Gamma }v_\gamma w_\gamma .
	\end{equation}
For $1 < p \le \infty$ this pairing identifies $L^p (\Gamma)$ with the dual $L^q (\Gamma)'$ as Banach spaces and $L^p (\Gamma)$ acquires a weak $*$-topology. In this topology $(v_n)$ converges to $v$ if and only if $\lim_{n\to \infty} \langle v_n , w \rangle = \langle v, w\rangle$ for all $w$ in $L^q (\Gamma)$.

Note that
	\begin{equation}
	\label{eq:adjoint}
	\begin{gathered}
\langle v,\rho_hw\rangle =\langle v,w\cdot h^*\rangle =\langle v\cdot h,w\rangle ,
	\\
\langle v,\lambda_hw\rangle =\langle v,h\cdot w\rangle =\langle h^*\cdot v,w\rangle
	\end{gathered}
	\end{equation}
for every $h\in L^1(\Gamma )$. Finally we put, for every $h\in L^1(\Gamma )$ and $1\le p\le \infty $,
	\begin{equation}
	\label{eq:notation}
K_p(h)=\{g\in L^p(\Gamma ):\rho_hg=0\},\enspace \enspace V_p(h)=\rho_h(L^p(\Gamma )).
	\end{equation}

\section{Algebraic $\Gamma $-actions}
	\label{s:actions}

	\begin{defi}
	\label{d:algebraic}
Let $\Gamma $ be a countable discrete group. An \emph{algebraic $\Gamma $-action} is a homomorphism $\alpha \colon \gamma \mapsto \alpha ^\gamma $ from $\Gamma $ into the group $\textup{Aut}(X)$ of continuous automorphisms of a compact abelian group $X$.

If $\alpha $ is an algebraic $\Gamma $-action on a compact abelian group $X$ then $\alpha $ is \emph{expansive} if there exists an open neighbourhood $\mathcal{U}$ of the identity in $X$ with $\bigcap_{\gamma \in\Gamma }\alpha ^\gamma (\mathcal{U})=0$.
	\end{defi}

In this section we consider algebraic $\Gamma $-actions of a very special nature. Fix a countable discrete group $\Gamma $ and consider the compact abelian group $X=\mathbb{T}^\Gamma $ consisting of all maps from $\Gamma $ to $\mathbb{T}=\mathbb{R}/\mathbb{Z}$ under point-wise addition. Each $x\in X$ is written as $(x_\gamma )=(x_\gamma ,\,\gamma \in\Gamma )$, where $x_\gamma \in\mathbb{T}$ denotes the value of $x$ at $\gamma \in\Gamma $. We identify the Pontryagin dual $\widehat{X}$ of $X$ with the group ring $\mathbb{Z}\Gamma$ under the pairing
	\begin{equation}
	\label{eq:pairing}
\langle f,x\rangle =e^{2\pi i\sum_{\gamma \in\Gamma }f_\gamma x_\gamma },
	\end{equation}
where $f=\sum_{\gamma \in\Gamma }f_\gamma \gamma \in\mathbb{Z}\Gamma$ and $x=(x_\gamma )\in X$.

The left and right shift actions $\lambda $ and $\rho $ of $\Gamma $ on $X$ are defined by
	\begin{equation}
	\label{eq:shift3}
(\lambda ^\gamma x)_{\gamma '}=x_{\gamma ^{-1}\gamma '}\enspace \enspace \textup{and}\enspace \enspace (\rho ^\gamma x)_{\gamma '}=x_{\gamma '\gamma }
	\end{equation}
for every $\gamma ,\gamma '\in\Gamma $ and $x\in X$. As in \eqref{eq:action} we can extend these actions of $\Gamma $ to commuting (left) actions $\lambda $ and $\rho $ of $\mathbb{Z}\Gamma$ on $X$ by setting
	\begin{equation}
	\label{eq:action2}
\lambda _fx=\sum_{\gamma \in\Gamma }f_\gamma \lambda ^\gamma x\enspace \enspace \textup{and}\enspace \enspace \rho _fx=\sum_{\gamma \in\Gamma }f_\gamma \rho ^{\gamma }x
	\end{equation}
for $f\in\mathbb{Z}\Gamma$ and $x\in X$ (cf. \eqref{eq:involution}). For every $f\in\mathbb{Z}\Gamma$, the maps $\lambda _f,\rho _f\colon X\linebreak[0]\longrightarrow X$ are continuous group homomorphisms which are dual to left and right multiplication by $f^*$ and $f$, respectively, on $\widehat{X}=\mathbb{Z}\Gamma$.

For the following discussion we set
	\begin{equation}
	\label{eq:LZ}
L^\infty (\Gamma ,\mathbb{Z})=\{w=(w_\gamma )\in L^\infty (\Gamma ):w_\gamma \in\mathbb{Z}\enspace \textup{for every}\enspace \gamma \in\Gamma \}.
	\end{equation}
We fix $f\in\mathbb{Z}\Gamma$, set
	\begin{equation}
	\label{eq:Xf}
X_f=\ker(\rho _f)=\{x\in X:\rho _fx=0\}= \widehat{\Z \Gamma / \Z \Gamma f},
	\end{equation}
and denote by
	\begin{equation}
	\label{eq:alphaf}
\alpha _f=\lambda |_{X_f}
	\end{equation}
the restriction to $X_f$ of the $\Gamma $-action $\lambda $ on $X$. Since $\alpha _f^\gamma \in\textup{Aut}(X_f)$ for every $\gamma \in\Gamma $, $\alpha _f$ is an algebraic $\Gamma $-action on the compact abelian group $X_f$.

\medskip The map $\eta \colon L^\infty (\Gamma )\longrightarrow X$, given by
	\begin{equation}
	\label{eq:eta}
\eta (w)_\gamma =w_\gamma \;(\textup{mod}\;1)
	\end{equation}
for every $w=(w_\gamma )\in L^\infty (\Gamma )$ and $\gamma \in\Gamma $, is a continuous surjective group homomorphism with
	\begin{equation}
	\label{eq:equivariant}
\eta \circ \lambda^\gamma =\lambda ^\gamma \circ \eta ,\enspace \enspace \eta \circ \rho^\gamma =\rho ^\gamma \circ \eta
	\end{equation}
for every $\gamma \in \Gamma $, and the \emph{linearization}
	\begin{equation}
	\label{eq:Wf}
W_f=\eta ^{-1}(X_f)=\rho_f^{-1}(L^\infty (\Gamma ,\mathbb{Z}))=\{w\in L^\infty (\Gamma ):\rho_fw\in L^\infty (\Gamma ,\mathbb{Z})\}
	\end{equation}
of $X_f$ is a weak$^*$-closed and $\lambda$-invariant subgroup with $\ker (\eta )=L^\infty (\Gamma ,\mathbb{Z})\subset W_f$.

	\begin{theo}
	\label{t:expansive}
Let $\Gamma $ be a countable group, $f\in\mathbb{Z}\Gamma$, and let $\alpha _f$ be the algebraic $\Gamma $-action on $X_f$ defined in \eqref{eq:Xf}--\eqref{eq:alphaf}. The following conditions are equivalent.

	\begin{enumerate}
	\item
The action $\alpha _f$ is expansive;
	\item
$K_\infty (f)=\{0\}$ \textup{(}cf. \eqref{eq:action} and \eqref{eq:notation}\textup{)};
	\item
$f$ is invertible in $L^1(\Gamma )$.
	\end{enumerate}
	\end{theo}

	\begin{proof}
This follows by combining \cite[Theorem 8.1]{ER} with \cite[Theorem 1]{Choi}. For the convenience of the reader we include the argument.

Let $d$ be the usual metric
	\begin{equation}
	\label{eq:d}
d(s_1,s_2)=\min\,\{|\tilde{s}_1-\tilde{s}_2|: \tilde{s}_i\in\mathbb{R}\enspace \textup{and}\enspace s_i=\tilde{s}_i\enspace (\textup{mod}\,1)\enspace \textup{for}\enspace i=1,2\}
	\end{equation}
on $\mathbb{T}$. If there exists a nonzero element $v\in K_\infty (f)$ then $\eta (cv)\in X_f$ for every $c\in \mathbb{R}$, and by choosing $|c|$ sufficiently small we can find, for every $\varepsilon >0$, a nonzero element $x^{(\varepsilon )}\in X_f$ with $d(x^{(\varepsilon )}_\gamma ,0)<\varepsilon $ for every $\gamma \in\Gamma $. This proves that $\alpha _f$ is nonexpansive.

Conversely, if $\alpha _f$ is nonexpansive, then we can find a nonzero element $x\in X_f$ with $d(x_\gamma ,0)<(3\|f\|_1)^{-1}$ for every $\gamma \in\Gamma $. We choose $\tilde{x}\in \eta ^{-1}(\{x\})\subset W_f=\eta ^{-1}(X_f)$ with $|\tilde{x}_\gamma |<(3\|f\|_1)^{-1}$ for every $\gamma \in\Gamma $. By definition of $X_f$, $\rho_f\tilde{x}\in L^\infty (\Gamma ,\mathbb{Z})$, and the smallness of the coordinates of $\tilde{x}$ implies that $\tilde{x}\in K_\infty (f)$. This proves that (1) $\Leftrightarrow$ (2).

If $K_\infty (f)=\{0\}$ then $V_1(f^*)$ is dense in $L^1(\Gamma )$ by \eqref{eq:adjoint} and the Hahn-Banach theorem. Since the group of units in $L^1(\Gamma )$ is open, $V_1(f^*)$ contains a unit. Hence there exists a $g\in L^1(\Gamma )$ with $g\cdot f=1$. By \cite[p. 122]{Kap}, $f$ is invertible in $L^1(\Gamma )$, which proves (3). The implication (3) $\Rightarrow$ (2) is obvious.
	\end{proof}

\section{Homoclinic points\label{s:homoclinic}}

	\begin{defi}
	\label{d:homoclinic}
Let $\alpha $ be an algebraic action of a countable discrete group $\Gamma $ on a compact abelian group $X$ with identity element $0$. A point $x\in X$ is \emph{$\alpha $-homoclinic} (or simply \emph{homoclinic}) if $\lim_{\gamma \to\infty }\alpha ^\gamma x=0$, i.e. if for every neighbourhood $U$ of $0$ in $X$ there is a finite subset $F$ of $\Gamma$ with $\alpha^{\gamma} x \in U$ for all $\gamma \in \Gamma \smallsetminus F$.

The set $\Delta _\alpha (X)$ of all $\alpha $-homoclinic points in $X$ is a subgroup of $X$, called the \emph{homoclinic group of $\alpha $}.

Following \cite{LS} we call an $\alpha $-homoclinic point $x\in X$ \emph{fundamental} if the homoclinic group $\Delta _\alpha (X)$ is generated by the orbit $\{\alpha ^\gamma x:\gamma \in\Gamma \}$ of $x$.
	\end{defi}

	\begin{prop}
	\label{p:homoclinic}
Let $\Gamma $ be a countable group and $f\in\mathbb{Z}\Gamma$ an element which is invertible in $L^1(\Gamma )$.

If $w_f^\Delta =f^{-1}\in L^1(\Gamma )$ and $\xi =\eta \circ \rho_{w_f^\Delta }\colon L^\infty (\Gamma ,\mathbb{Z})\longrightarrow X_f$ \textup{(}cf. \eqref{eq:action} and \eqref{eq:Xf}--\eqref{eq:eta}\textup{)}, then $\xi $ is a surjective group homomorphism with the following properties.
	\begin{enumerate}
	\item
$\ker(\xi )=\rho_f(L^\infty (\Gamma ,\mathbb{Z}))$;
	\item
$\xi \circ \lambda^\gamma =\alpha _f^\gamma \circ \xi $ for every $\gamma \in\Gamma $;
	\item
$\xi $ is continuous in the weak$^*$-topology on closed, bounded subsets of $L^\infty (\Gamma ,\mathbb{Z})$.
	\end{enumerate}
	\end{prop}

	\begin{proof}
We set $\tilde{w}=(f^*)^{-1}=(f^{-1})^*$. By definition, $\rho_f\tilde{w}=\tilde{w}\cdot f^*=e(1)$ and hence $\tilde{w}\in W_f$ and $x_f^\Delta =\eta (\tilde{w})\in X_f$ (cf. \eqref{eq:Wf}). Since $\tilde{w}\in L^1(\Gamma )$, $x_f^\Delta \in\Delta _{\alpha _f}(X_f)$.

The $\lambda$-invariance of $W_f$ implies that $\rho_{w_f^\Delta }h=\lambda_h\tilde{w}\in W_f$ for every $h\in \mathbb{Z}\Gamma$. Since $W_f$ is weak$^*$-closed and $\rho_{w_f^\Delta }$ is weak$^*$-continuous on bounded subsets of $L^\infty (\Gamma ,\mathbb{Z})$, it follows that
	$$
\rho_{w_f^\Delta }(L^\infty (\Gamma ,\mathbb{Z}))\subset W_f.
	$$

In order to prove that $\rho_{w_f^\Delta }(L^\infty (\Gamma ,\mathbb{Z}))=W_f$ we fix $w\in W_f$, set $v=\rho_fw\in L^\infty (\Gamma ,\mathbb{Z})$, and obtain that $w=\rho_{w_f^\Delta }v$.

The group homomorphism
	\begin{equation}
	\label{eq:xi}
\xi =\eta \circ \rho_{w_f^\Delta }\colon L^\infty (\Gamma ,\mathbb{Z})\longrightarrow X_f
	\end{equation}
is thus surjective, and the equivariance of $\xi $ is obvious.

If $B\subset L^\infty (\Gamma ,\mathbb{Z})$ is a closed, bounded subset, then the weak$^*$-topology coincides with the topology of coordinate-wise convergence, and $\xi $ is obviously continuous in that topology.
	\end{proof}

	\begin{theo}
	\label{t:homoclinic}
Let $\Gamma $ be a countable residually finite discrete group and $f\in\mathbb{Z}\Gamma$. If the algebraic $\Gamma $-action $\alpha _f$ on the compact abelian group $X_f$ in \eqref{eq:Xf}--\eqref{eq:alphaf} is expansive, and if $w_f^\Delta =f^{-1}\in L^1(\Gamma )$ is the inverse of $f$ in $L^1(\Gamma )$ described in Theorem \ref{t:expansive}, then $x_f^\Delta =\eta \bigl((w_f^\Delta )^*\bigr)\in X_f$ is a fundamental homoclinic point of $\alpha _f$.
	\end{theo}

	\begin{proof}
Since $\tilde{w}=(w_f^\Delta )^*\in L^1(\Gamma )$, $x_f^\Delta \in\Delta _{\alpha _f}(X_f)$. If $x\in\Delta _{\alpha _f}(X_f)$ is an arbitrary $\alpha _f$-homoclinic point, then we can find a $w\in\eta ^{-1}(\{x\})$ with $\lim_{\gamma \to\infty }w_\gamma =0$ and hence with $v=\rho_f(w)\in L^1(\Gamma ,\mathbb{Z})=L^1(\Gamma )\cap L^\infty (\Gamma ,\mathbb{Z})$. The point $w=\rho_{w_f^\Delta }v$ lies in the group generated by the $\lambda$-orbit $\{\lambda^\gamma \tilde{w}:\gamma \in\Gamma \}$ of $\tilde{w}$, and by applying $\eta $ we see that $x$ lies in the group generated by the $\alpha _f$-orbit of $x_f^\Delta $. This proves that $x_f^\Delta $ is fundamental.
	\end{proof}

	\begin{theo}[Specification Theorem]
	\label{t:specification}
Let $\Gamma $ be a countable residually finite group and $f\in\mathbb{Z}\Gamma$ an element such that the algebraic $\Gamma $-action $\alpha _f$ on the compact abelian group $X_f$ in \eqref{eq:Xf}--\eqref{eq:alphaf} is expansive. Then there exists, for every $\varepsilon >0$, a finite subset $F_\varepsilon \subset \Gamma $ with the following property: if $C_1,C_2$ are subsets of $\Gamma $ with $F_\varepsilon C_1\cap F_\varepsilon C_2=\varnothing $, then we can find, for every pair of points $x^{(1)},x^{(2)}\in X_f$, a point $y\in X_f$ with $d(x_\gamma ^{(i)},y_\gamma )<\varepsilon $ for every $\gamma \in C_i,\,i=1,2$ \textup{(}cf. \eqref{eq:d}\textup{)}.
	\end{theo}

For the proof of Theorem \ref{t:specification} we need an elementary lemma.

	\begin{lemm}
	\label{l:l6}
For every $x\in X_f$ there exists an element $v\in L^\infty (\Gamma ,\mathbb{Z})$ with $\xi (v)=x$ and $\|v\|_\infty \le\|f\|_1/2$.
	\end{lemm}

	\begin{proof}
Choose $w\in W_f=\eta ^{-1}(X_f)\subset L^\infty (\Gamma )$ with $\eta (w)=x$ and $-1/2\le w_\gamma <1/2$ for every $\gamma \in\Gamma $. Then $v=\rho_fw\in L^\infty (\Gamma ,\mathbb{Z})$, $\|v\|_\infty \le \|f\|_1/2$ and $\xi (v)=x$ by Proposition \ref{p:homoclinic}.
	\end{proof}

	\begin{proof}[Proof of Theorem \ref{t:specification}]
Consider the point $w_f^\Delta =f^{-1}\in L^1(\Gamma )$ described in Theorem \ref{t:expansive}. We can find, for every $\varepsilon >0$, a finite set $F_\varepsilon \subset \Gamma $ with $\sum_{\gamma \in\Gamma \smallsetminus F_\varepsilon }|(w_f^\Delta )_\gamma |<\varepsilon /\|f\|_1$. Put
	$$
\tilde{w} _\gamma =
	\begin{cases}
(w_f^\Delta )_\gamma &\textup{if}\enspace \gamma \in F_\varepsilon ,
	\\
0&\textup{otherwise}
	\end{cases}
	$$
and note that $|(\rho_{w_f^\Delta }y)_\gamma -(\rho_{\tilde{w}}y)_\gamma |<\varepsilon /2$ for every $\gamma \in\Gamma $ and every $y\in L^\infty (\Gamma ,\mathbb{Z})$ with $\|y\|_\infty \le \|f\|_1/2$.

Lemma \ref{l:l6} allows us to find points $v^{(i)}\in L^\infty (\Gamma ,\mathbb{Z})$ such that $\xi (v^{(i)})=x^{(i)}$ and $\|v^{(i)}\|\le \|f\|_1/2$ for $i=1,2$. Let $v\in L^\infty (\Gamma ,\mathbb{Z})$ be a point with $\|v\|_\infty \le \|f\|_1/2$ and $v_\gamma =v_\gamma ^{(i)}$ for $\gamma \in C_iF_\varepsilon ^{-1},\,i=1,2$, where $F_\varepsilon ^{-1}=\{\gamma ^{-1}:\gamma \in F_\varepsilon \}$. Then $d(x_\gamma ^{(i)},\eta (\rho_{\tilde{w}}v^{(i)})_\gamma )<\varepsilon /2$ for every $\gamma \in\Gamma $ and $(\rho_{\tilde{w}}v^{(i)})_\gamma =(\rho_{\tilde{w}}v)_\gamma $ for every $\gamma \in C_i,\,i=1,2$. By setting $y=\xi (v)$ we obtain that $d(x_\gamma ^{(i)},y_\gamma )<\varepsilon $ for every $\gamma \in C_i,\,i=1,2$. This proves our claim.
	\end{proof}

	\begin{theo}
	\label{t:dense}
Let $\Gamma $ be a countable residually finite group and $f\in\mathbb{Z}\Gamma$ an element such that the algebraic $\Gamma $-action $\alpha _f$ on the compact abelian group $X_f$ in \eqref{eq:Xf}--\eqref{eq:alphaf} is expansive. Then the homoclinic group $\Delta _{\alpha _f}(X_f)$ is countable and dense in $X_f$.
	\end{theo}

	\begin{proof}
The countability of $\Delta _{\alpha _f}(X_f)$ is proved exactly as in \cite[Lemma 3.2]{LS}. In order to verify that $\Delta _{\alpha _f}(X_f)$ is dense in $X_f$ we consider the continuous surjective group homomorphism $\xi \colon L^\infty (\Gamma ,\mathbb{Z})\longrightarrow X_f$ in Proposition \ref{p:homoclinic}. We set $B=\{v\in L^\infty (\Gamma ,\mathbb{Z}):\|v\|_\infty \le\|f\|_1/2\}$ and note that $\xi (B)=X_f$ by Lemma \ref{l:l6}, and that the restriction of $\xi $ to $B$ is continuous in the weak$^*$-topology. Since $B\cap L^1(\Gamma ,\mathbb{Z})$ is countable and weak$^*$-dense in $B$, the countable set $\xi (B\cap L^1(\Gamma ,\mathbb{Z}))\subset \Delta _{\alpha _f}(X_f)$ is dense in $X_f$.
	\end{proof}

We end this section with a lemma which implies that the point $w_f^\Delta \in L^1(\Gamma )$ appearing in the Theorems \ref{t:homoclinic} and \ref{t:specification} decays rapidly.

	\begin{prop}
	\label{p:expansive}
Let $\Gamma $ be a countable discrete group and $f\in\mathbb{Z}\Gamma$. If $\alpha _f$ is expansive there exists, for every $L>0$ and $\varepsilon >0$, a finite subset $F(L,\varepsilon )\subset \Gamma $ with the following property: if $w\in W_f$ satisfies that $\|w\|_\infty \le L$ and $(\rho_fw)_\gamma =0$ for every $\gamma \in F(L,\varepsilon )$, then $|w_1|<\varepsilon $.
	\\
If the group $\Gamma $ is residually finite, finitely generated and has polynomial growth, then the point $w_f^\Delta $ has exponential decay in the word metric on $\Gamma $.
	\end{prop}

	\begin{proof}
We argue by contradiction and assume that there exist an $\varepsilon >0$, an increasing sequence $(F_n)$ of finite subsets in $\Gamma $ with $\bigcup_{n\ge1}F_n=\Gamma $ and a sequence $(w^{(n)},\,n\ge1)$ in $W_f$ such that, for every $n\ge1$, $\|w^{(n)}\|_\infty \le L$, $|w^{(n)}_1|>\varepsilon $ and $\rho_fw^{(n)}_\gamma =0$ for every $\gamma \in F_n$. If $w^{(0)}$ is the limit of a weak$^*$-convergent subsequence of $(w^{(n)})$, then $|w^{(0)}_1|\ge\varepsilon $ and $w^{(0)}\in\ker(\rho_f)$, contrary to our hypothesis that $\alpha _f$ is expansive and $\ker(\rho_f)$ is therefore equal to $0$.

Now suppose that $\Gamma $ is residually finite, finitely generated and has polynomial growth. We choose and fix a finite symmetric set of generators $F$ of $\Gamma $ and write $\delta _F$ for the word-metric on $\Gamma $. For every $\gamma \in\Gamma $ we denote by $\ell _F(\gamma )$ the length of the shortest expression of $\gamma \in\Gamma $ in terms of elements of $F$, where $\ell _F(1)=0$. Since $\Gamma $ has polynomial growth there exist constants $c,M\ge1$ such that the set
	\begin{equation}
	\label{eq:BF}
B_F(n)=\{\gamma \in\Gamma :\ell _F(\gamma )\le n\}
	\end{equation}
has cardinality $\le cM^n$ for every $n\ge0$.

The point $\tilde{w}=(w_f^\Delta )^*\in L^1(\Gamma )$ satisfies that $\rho_f\tilde{\omega }=e(1)$. We set $L=\|\tilde{w}\|_\infty $ and use the first part of this proposition to find an $R\ge1$ such that $|w_1|<L/2$ for every $w\in W_f$ with $\|w\|_\infty \le L$ and $\rho_fw\equiv0$ on $B_F(R)$. Our choice of $R$ guarantees that $|\tilde{w}_\gamma |<L/2$ for every $\gamma \in B_F(2R)\smallsetminus B_F(R)$ and, by induction, that $|\tilde{w}_\gamma |<L/2^k$ for every $k\ge1$ and every $\gamma \in B_F(kR)\smallsetminus B_F((k-1)R)$.

This proves exponential decay of the coordinates of $\tilde{w}$ (and hence of $w_f^\Delta $) in the word metric on $\Gamma $.
	\end{proof}

\section{Entropy and periodic points}
	\label{s:entropy}

Let $\Gamma $ be a countable discrete group and $K\subset \Gamma $ a finite set. A finite set $Q\subset \Gamma $ is \emph{left $(K,\varepsilon )$-invariant} if
	$$
\sum_{\gamma \in K}|\gamma Q\triangle Q|/|Q|<\varepsilon ,
	$$
and \emph{right $(K,\varepsilon )$-invariant} if
	$$
\sum_{\gamma \in K}|Q\gamma \triangle Q|/|Q|<\varepsilon .
	$$
If $Q$ satisfies both these conditions it is \emph{$(K,\varepsilon )$-invariant}.

A sequence $(Q_n,\,n\ge1)$ of finite subsets of $\Gamma $ is a \emph{left F{\o}lner sequence} if there exists, for every finite subset $K\subset \Gamma $ and every $\varepsilon >0$, an $N\ge1$ such that $Q_n$ is left $(K,\varepsilon )$-invariant for every $n\ge N$. The definitions of \emph{right} and \emph{two-sided F{\o}lner sequences} are analogous. The group $\Gamma $ is \emph{amenable} if it has a left F{\o}lner sequence. If $\Gamma $ is amenable it also has right and two-sided F{\o}lner sequences.

Let $\Gamma$ be a countable residually finite discrete group. If $(\Gamma _n,\,n\ge1)$ is a sequence of finite index normal subgroups in $\Gamma $ we say that
	\begin{equation}
	\label{eq:limit}
\lim_{n\to\infty }\Gamma _n=\{1\}
	\end{equation}
if we can find, for every finite set $K\subset \Gamma $, an $N\ge1$ with $\Gamma _n\cap(K^{-1}K)=\{1\}$ for every $n\ge N$. Clearly, such sequences exist.

	\begin{theo}
	\label{t:positive}
Let $\Gamma $ be a countable residually finite amenable group, $f\in\mathbb{Z}\Gamma$, and let $\alpha _f$ be the $\Gamma $-action on $X_f$ defined in \eqref{eq:Xf}--\eqref{eq:alphaf}. If $X_f\ne0$ and $\alpha _f$ is expansive then $h(\alpha _f)>0$.
	\end{theo}

	\begin{proof}
Let $(\Gamma _n,\,n\ge1)$ be a decreasing sequence of finite index normal subgroups of $\Gamma $ with $\bigcap_{n\ge1}\Gamma _n=\{1\}$.

For notational simplicity set $\tilde{w}=(w_f^\Delta )^*$ (cf. Theorem \ref{t:homoclinic}). The homoclinic group $\Delta _{\alpha _f}(X_f)$ is dense in $X_f$ by Theorem \ref{t:dense} and, since $X_f\ne0$ by assumption, the fundamental homoclinic point $\tilde{x}=x_f^\Delta =\eta (\tilde{w})$ is nonzero. Hence there exist a $\gamma _0\in \Gamma $ and a finite subset $F\subset \Gamma $ with $\tilde{w}_{\gamma _0}\notin\mathbb{Z}$ and $\sum_{\gamma \in\Gamma \smallsetminus F}|\tilde{w}_\gamma |<d(\tilde{x}_{\gamma _0},0)/2$, where $d(\cdot ,\cdot )$ is the metric on $\mathbb{T}$ defined in \eqref{eq:d}). We choose $n$ sufficiently large so that the sets $\gamma F,\,\gamma \in\Gamma _n$, are all disjoint. For every $\omega \in\Omega =\{0,1\}^{\Gamma _n}$ we define a point $x^{(\omega )}\in X_f$ by setting $w^{(\omega )}=\sum_{\gamma \in\Gamma _n}\omega _\gamma \lambda^\gamma \tilde{w}$ and $x^{(\omega )}=\eta (w^{(\omega )})=\sum_{\gamma \in\Gamma _n}\omega _\gamma \lambda ^\gamma \tilde{x}$.

We fix $\gamma \in\Gamma _n$ for the moment and consider two points $\omega ,\omega '\in\Omega $ with $\omega _\gamma \ne \omega _\gamma '$. Then $d(x^{(\omega )}_{\gamma \gamma _0},x^{(\omega ')}_{\gamma \gamma _0})>d(\tilde{x}_{\gamma _0},0)/2$.

If $(Q_m,\,m\ge1)$ is a left F{\o}lner sequence in $\Gamma $ then
	$$
\lim_{m\to\infty }\frac{|Q_m\cap \Gamma _n\gamma _0|}{|Q_m|}=|\Gamma /\Gamma _n|^{-1}.
	$$
We conclude that the set $\{x^{(\omega )}:\omega \in\Omega \}$ contains, for every $m\ge1$, a $(Q_m,d(\tilde{x}_{\gamma _0},0)/2)$-separated set of cardinality $2^{|Q_m\cap \Gamma _n\gamma _0|}$, and hence that $h(\alpha _f)\linebreak[0]\ge |\Gamma /\Gamma _n|^{-1}\cdot \log 2>0$.
	\end{proof}

In order to find out more about the actual value of $h(\alpha _f)$ we take a look at the periodic points of $\alpha _f$. Fix a subgroup $\Gamma '\subset \Gamma $ and denote by
	\begin{equation}
	\label{eq:Fix}
\textup{Fix}_{\Gamma '}(X_f)=\{x\in X_f:\alpha _f^\gamma x=x\enspace \textup{for every}\enspace \gamma \in\Gamma '\}
	\end{equation}
the subgroup of \emph{$\Gamma '$-invariant points} in $X_f$. Clearly, $\textup{Fix}_{\Gamma '}(X_f)$ is $\Gamma '$-invariant, and $\textup{Fix}_{\Gamma '}(X_f)$ is $\Gamma $-invariant if and only if $\Gamma '$ is a normal subgroup of $\Gamma $. Next we set
	\begin{equation}
	\label{eq:Fix2}
	\begin{gathered}
L^\infty (\Gamma )^{\Gamma '}=\{w\in L^\infty (\Gamma ):\lambda^\gamma w=w\enspace \textup{for every}\enspace \gamma \in \Gamma '\},
	\\
W_f^{\Gamma '}=W_f\cap L^\infty (\Gamma )^{\Gamma '},
	\\
L^\infty (\Gamma ,\mathbb{Z})^{\Gamma '}=L^\infty (\Gamma ,\mathbb{Z})\cap L^\infty (\Gamma )^{\Gamma '}.
	\end{gathered}
	\end{equation}
We write $\xi \colon L^\infty (\Gamma ,\mathbb{Z})\longrightarrow X_f$ for the group homomorphism described in Proposition \ref{p:homoclinic}.

	\begin{prop}
	\label{p:periodic}
Let $\Gamma $ be a countable residually finite group, $f\in\mathbb{Z}\Gamma$, and let $\alpha _f$ be the $\Gamma $-action on $X_f$ defined in \eqref{eq:Xf}--\eqref{eq:alphaf}. For every subgroup $\Gamma '\subset \Gamma $ of finite index,
	\begin{equation}
	\label{eq:Fix3}
\textup{Fix}_{\Gamma '}(X_f)=\xi (L^\infty (\Gamma ,\mathbb{Z})^{\Gamma '}) \cong L^\infty (\Gamma ,\mathbb{Z})^{\Gamma '}/\rho_f(L^\infty (\Gamma ,\mathbb{Z})^{\Gamma '}).
	\end{equation}
	\end{prop}

	\begin{proof}
From the equivariance of $\xi $ it is clear that that $\xi (L^\infty (\Gamma ,\mathbb{Z})^{\Gamma '})\subset \textup{Fix}_{\Gamma '}(X_f)$. Conversely, if $x\in \textup{Fix}_{\Gamma '}(X_f)$, then we can find $w\in W_f^{\Gamma '}\subset L^\infty (\Gamma )^{\Gamma '}$ with $\eta (w)=x$ (cf. \eqref{eq:eta}), and the point $v=\rho_fw\in L^\infty (\Gamma ,\mathbb{Z})^{\Gamma '}$ satisfies that $\xi (v)=x$. This proves that $\xi (L^\infty (\Gamma ,\mathbb{Z})^{\Gamma '})=\textup{Fix}_{\Gamma '}(X_f)$, and Proposition \ref{p:homoclinic} guarantees that $\ker(\xi )\cap L^\infty (\Gamma ,\mathbb{Z})^{\Gamma '}=\rho_f(L^\infty (\Gamma ,\mathbb{Z})^{\Gamma '})$. It follows that $\textup{Fix}_{\Gamma '}(X_f)=\xi (L^\infty (\Gamma ,\mathbb{Z})^{\Gamma '})\cong L^\infty (\Gamma ,\mathbb{Z})^{\Gamma '}/\rho_f(L^\infty (\Gamma ,\mathbb{Z})^{\Gamma '})$, as claimed.
	\end{proof}

	\begin{coro}
	\label{c:periodic}
If $|\Gamma '\backslash \Gamma |<\infty $ then $|\textup{Fix}_{\Gamma '}(X_f)|=|\negthinspace\det(\rho_f|_{L^\infty (\Gamma )^{\Gamma '}})|$.
	\end{coro}

	\begin{proof}
If the coset space $\Gamma '\backslash \Gamma $ is finite then $L^\infty (\Gamma )^{\Gamma '}\cong \mathbb{R}^{\Gamma '\backslash \Gamma }$ and we denote by $\rho_f|_{L^\infty (\Gamma )^{\Gamma '}}$ the restriction to $L^\infty (\Gamma )^{\Gamma '}\cong\mathbb{Z}^{\Gamma '\backslash \Gamma }$ of the linear map $\rho_f\colon L^\infty (\Gamma )\longrightarrow L^\infty (\Gamma )$. Then $\rho_f(L^\infty (\Gamma ,\mathbb{Z})^{\Gamma '})\subset L^\infty (\Gamma ,\mathbb{Z})^{\Gamma '}$ and the absolute value of the determinant $\negthinspace\det(\rho_f|_{L^\infty (\Gamma )^{\Gamma '}})$ is equal to
	\\
$|L^\infty (X,\mathbb{Z})^{\Gamma '}\slash\linebreak[0]\rho_f(L^\infty (\Gamma ,\mathbb{Z})^{\Gamma '})|=|X_f^{\Gamma '}|$.
	\end{proof}

	\begin{rema}
	\label{r:periodic}
Since $\alpha _f$ is expansive, $\ker(\rho_f)=0$ by Theorem \ref{t:expansive}, and the proof of Theorem \ref{t:expansive} shows that there exists, for every pair of distinct points $x,x'\in X_f$, a $\gamma \in\Gamma $ with $d(x_\gamma ,x_\gamma ')\ge(3\|f\|_1)^{-1}$ (cf. \eqref{eq:d}).

If $\Gamma '\subset \Gamma $ is a subgroup with finite index, and if $Q\subset \Gamma $ is a \emph{fundamental domain} of the right coset space $\Gamma '\backslash \Gamma $, i.e. a finite subset such that $\{\gamma Q\colon \gamma \in \Gamma '\}$ is a partition of $\Gamma $, then the preceding paragraph implies that $\textup{Fix}_{\Gamma '}(X_f)$ is \emph{$(Q,(3\|f\|_1)^{-1})$-separated} in the sense that there exists, for any two distinct points $x,x'\in X_f^{\Gamma '}$, a $\gamma \in Q$ with $d(x_\gamma ,x_\gamma ')\ge(3\|f\|_1)^{-1}$.
	\end{rema}

For the terminology used in our next proposition we refer to the beginning of this section and to Remark \ref{r:periodic}.

	\begin{prop}[\cite{Weiss}]
	\label{p:Weiss}
Let $\Gamma $ be a countable residually finite amenable group and let $(\Gamma _n,\,n\ge1)$ be a sequence of finite index normal subgroups with $\lim_{n\to\infty }\Gamma _n=\{1\}$ \textup{(}cf. \eqref{eq:limit}\textup{)}. Then there exists, for every finite subset $K\subset \Gamma $ and every $\bar{\varepsilon }>0$, an integer $M=M(K,\bar{\varepsilon })\ge1$ such that every $\Gamma _n$ with $n\ge M$ has a $(K,\bar{\varepsilon })$-invariant fundamental domain $Q_n$ of the coset space $\Gamma /\Gamma _n$.
	\end{prop}

	\begin{proof}
This is a slight reformulation of \cite[Theorem 1]{Weiss} requiring only very minor changes in the proof. We describe these changes briefly, using essentially the same notation and terminology as in \cite{Weiss}.

Fix a finite set $K\subset \Gamma $ and $\bar{\varepsilon }>0$ and let $\varepsilon $ and $\eta $ be sufficiently small positive numbers whose sizes will be clear by examining the course the proof. Choose $N\ge1$ with
	$$
\bigl(1-\tfrac{\varepsilon }2\bigr)^N<\tfrac{\bar{\varepsilon }}{200}
	$$
and use the amenability of $\Gamma $ to find an increasing sequence $(F_j,\,j\ge1)$ such that
	\begin{enumerate}
	\item[(i)]
$F_1$ is $(K,\varepsilon )$-invariant and $1\in F_1$,
	\item[(ii)]
for $1\le j <N$, $F_{j+1}\supset F_j$ and $F_{j+1}$ is $(F_jF_j^{-1},\eta )$-invariant.
	\end{enumerate}
Since $\bigcap_{n\ge1}\Gamma _n=\{1\}$ we can find an $M\ge1$ with $\Gamma _M\cap F_N^{-1}F_N=\{1\}$. We fix $n\ge M$ and write $\theta \colon \Gamma \longrightarrow \Gamma /\Gamma _n=G$ for the quotient map which is injective on $F_N$.

The argument in the proof of \cite[Theorem 1]{Weiss} allows us to find --- for appropriately chosen $\varepsilon ,\eta $ --- finite subsets $\{\gamma _{i,j}:i=1,\dots ,m_j,\,j=1,\dots ,N\}$ and $\{\gamma _{i,j}':i=1,\dots ,m_j',\,j=1,\dots ,N\}$ of $\Gamma $ and, for each $j=1,\dots ,N$, sets $F_{i,j}\subset F_j,\,i=1,\dots ,m_j$, $F_{i,j}'\subset F_j,\,i=1,\dots ,m_j'$, such that the following conditions are satisfied.
	\begin{enumerate}
	\item
For every $i,j$, $|F_{i,j}|\ge(1-\varepsilon )|F_j|$ and $|F_{i,j}'|\ge(1-\varepsilon )|F_j|$;
	\item
For every $i,j$, the sets $F_{i,j}$ and $F_{i,j}'$ are $(K,\bar{\varepsilon }/2)$-invariant;
	\item
The sets $\theta (F_{i,j}\gamma _{i,j}),\,i=1,\dots ,m_j,\,j=1,\dots ,N$, are disjoint;
	\item
The sets $\theta (\gamma _{i,j}'F_{i,j}'),\,i=1,\dots ,m_j',\,j=1,\dots ,N$, are disjoint;
	\item
The sets
	$$
E=\bigcup_{j=1}^N\bigcup_{i=1}^{m_j}\theta (F_{i,j}\gamma _{i,j}),\enspace E'=\bigcup_{j=1}^N\bigcup_{i=1}^{m_j'}\theta (\gamma _{i,j}'F_{i,j}')
	$$
satisfy that
	$$
|E|\ge (1-\tfrac{\bar{\varepsilon }}{100})|G|,\enspace |E'|\ge (1-\tfrac{\bar{\varepsilon }}{100})|G|.
	$$
	\end{enumerate}
We set $F=E\cap E'$ and choose a set $Q\supset F$ such that $|Q|=|G|$ and $\theta (Q)=G$. Then $Q$ is a $(K,\bar{\varepsilon })$-invariant fundamental domain of the coset space $\Gamma /\Gamma _n$.
	\end{proof}

	\begin{coro}
	\label{c:Weiss}
Let $\Gamma $ be a countable residually finite amenable group and let $(\Gamma _n,\,n\ge1)$ be a sequence of finite index normal subgroups with $\lim_{n\to\infty }\Gamma _n=\{1\}$. Then there exists a F{\o}lner sequence $(Q_n,\,n\ge1)$ such that $Q_n$ is a fundamental domain of $\Gamma /\Gamma _n$ for every $n\ge1$.
	\end{coro}

	\begin{theo}
	\label{t:entropy}
Let $\Gamma $ be a countable residually finite amenable group and let $(\Gamma _n,\,n\ge1)$ be a sequence of finite index normal subgroups with $\lim_{n\to\infty }\Gamma _n=\{1\}$.

If $f\in\mathbb{Z}\Gamma$, and if the algebraic $\Gamma $-action $\alpha _f$ on $X_f$ in \eqref{eq:Xf}--\eqref{eq:alphaf} is expansive, then
	\begin{equation}
	\label{eq:entropy}
	\begin{aligned}
h(\alpha _f)&=\lim_{n\to\infty }\frac{1}{|\Gamma /\Gamma _n|}\log\,|\textup{Fix}_{\Gamma _n}(X_f)|
	\\
&=\lim_{n\to\infty }\frac{1}{|\Gamma /\Gamma _n|}\log\,|L^\infty (\Gamma ,\mathbb{Z})^{\Gamma _n}/\rho_f(L^\infty (\Gamma ,\mathbb{Z})^{\Gamma _n})|
	\\
&=\lim_{n\to\infty }\frac{1}{|\Gamma /\Gamma _n|}\log\,|\negthinspace\det (\rho_f|_{L^\infty (\Gamma )^{\Gamma _n}})|,
	\end{aligned}
	\end{equation}
where $\textup{Fix}_{\Gamma _n}(X_f)\subset X_f$ is the subgroup of $\Gamma _n$-periodic points in \eqref{eq:Fix} and $L^\infty (\Gamma )^{\Gamma _n}\linebreak[0]\subset L^\infty (\Gamma )$ is defined in \eqref{eq:Fix2} \textup{(}cf. Corollary \ref{c:periodic}\textup{)}.
	\end{theo}

	\begin{proof}
We choose a F{\o}lner sequence $(Q_n,\,n\ge1)$ in $\Gamma $ such that $Q_n$ is a fundamental domain of $\Gamma /\Gamma _n$ for every $n\ge1$ (cf. Corollary \ref{c:Weiss}). Proposition \ref{p:periodic} and Corollary \ref{c:periodic} show that there exists, for every $n\ge1$, a $(Q_n,(3\|f\|_1)^{-1})$-separated set (in the sense of Remark \ref{r:periodic}) of cardinality
	$$
|\textup{Fix}_{\Gamma _n}(X_f)|=|L^\infty (\Gamma ,\mathbb{Z})^{\Gamma _n}/\rho_f(L^\infty (\Gamma ,\mathbb{Z})^{\Gamma _n})|=|\negthinspace\det (\rho_f|_{L^\infty (\Gamma )^{\Gamma _n}})|.
	$$
Since $(Q_n,\,n\ge1)$ is F{\o}lner and $|Q_n|=|\Gamma /\Gamma _n|$ this implies that
	$$
h(\alpha _f)\ge\limsup_{n\to\infty }\frac{1}{|Q_n|}\log\,|\textup{Fix}_{\Gamma _n}(X_f)|.
	$$

Conversely, let $\delta >0$, $\varepsilon <\delta /3$, and let $F_\varepsilon $ be a finite symmetric set with $\sum_{\gamma \in\Gamma \smallsetminus F_\varepsilon }|w_\gamma ^\Delta |\linebreak[0]<\varepsilon /\|f\|_1$ (cf. the proof of Theorem \ref{t:specification}). The sets $P_n=Q_n\cap\,\bigcap_{\gamma \in F_\varepsilon }Q_n\gamma ,\, n\ge1$, form a F{\o}lner sequence with $\lim_{n\to\infty }\frac{|P_n|}{|Q_n|}=1$.

We fix $n\ge1$ for the moment and choose a maximal set $S_{n,\delta }\subset X_f$ which is $(P_n,\delta )$-separated in the sense of Remark \ref{r:periodic}. For every $x\in S_{n,\delta }$ we find $w(x)\in W_f\subset L^\infty (\Gamma )$ with $\|w(x)\|_\infty \le\|f\|_1/2$ and $\eta (w(x))=x$ (Lemma \ref{l:l6}) and write $v(x)\in L^\infty (\Gamma ,\mathbb{Z})^{\Gamma _n}$ for the unique point with $v(x)_\gamma =(\rho_fw(x))_\gamma $ for every $\gamma \in Q_n$. Our choice of $F_\varepsilon $ implies that the points $\{\xi (v(x)):x\in S_{n,\delta }\}\subset \textup{Fix}_{\Gamma _n}(X_f)$ are $(P_n,\delta /3)$-separated and therefore distinct, and Proposition \ref{p:periodic} shows that $|S_{n,\delta }|\le |\textup{Fix}_{\Gamma _n}(X_f)|$.

Since $(P_n,\,n\ge1)$ is F{\o}lner and $\lim_{n\to\infty }\frac{|P_n|}{|Q_n|}=1$ this implies that
	$$
h(\alpha _f)=\lim_{n\to\infty }\frac{1}{|P_n|}\log\,|S_{n,\delta }|\le \liminf_{n\to\infty }\frac{1}{|Q_n|}\log|\textup{Fix}_{\Gamma _n}(X_f)|,
	$$
which completes the proof of \eqref{eq:entropy}.
	\end{proof}

	\begin{coro}
	\label{c:f*}
Let $\Gamma $ be a countable residually finite amenable group, and let $f,g\in\mathbb{Z}\Gamma$ be elements such that the algebraic $\Gamma $-actions $\alpha _f$ and $\alpha _g$ on $X_f$ and $X_g$ in \eqref{eq:Xf}--\eqref{eq:alphaf} are expansive. If $\alpha _{f^*}$ and $\alpha _{f\cdot g}$ are given as in \eqref{eq:Xf}--\eqref{eq:alphaf} with $f$ replaced by $f^*$ and $f\cdot g$, respectively, then $\alpha _{f^*}$ and $\alpha _{f\cdot g}$ are expansive, $h(\alpha _{f^*})=h(\alpha _f)$ and $h(\alpha _{f\cdot g})=h(\alpha _f)+h(\alpha _g)$.
	\end{coro}

	\begin{proof}
The expansiveness of $\alpha _{f^*}$ and $\alpha _{f\cdot g}$ is clear from Theorem \ref{t:expansive}, and the entropy formulae follow from the usual properties of determinants (cf. \eqref{eq:entropy}).
	\end{proof}

\section{Entropy and Fuglede--Kadison determinants}
	\label{s:kadison}

In this section $L^p (\Gamma,\mathbb{C})$ will denote the \textit{complex} $L^p$-space of $\Gamma$ for $1 \le p \le \infty$ with its conjugate linear involution $w\mapsto w^*$ given by $(w^*)_{\gamma} = \overline{w}_{\gamma^{-1}},\;\gamma \in\Gamma $.

The von~Neumann algebra $\Nh\Gamma$ of a discrete group $\Gamma$ can be defined as the algebra of left $\Gamma$-equivariant bounded operators of $L^2 (\Gamma,\mathbb{C})$ to itself. Thus a bounded operator $A$ on $L^2 (\Gamma,\mathbb{C})$ belongs to $\Nh\Gamma$ if and only if we have \[ \lambda^{\gamma} A (v) = A (\lambda^{\gamma} (v)) \] for all $v$ in $L^2 (\Gamma,\mathbb{C})$ and all $\gamma$ in $\Gamma$.

The homomorphism of $\C$-algebras with involution: \[ \rho : L^1 (\Gamma ,\mathbb{C}) \longrightarrow \Nh\Gamma \] mapping $f$ to the operator $\rho_f$ with $\rho_f (v) = v \cdot f^*$ is injective because $\rho_f (e (1)) = f^*$. In the following we will sometimes view $\rho$ as an inclusion both of $L^1 (\Gamma ,\mathbb{C})$ and of $\C \Gamma$ into $\Nh\Gamma$ and omit it from the notation.

The von~Neumann trace on $\Nh\Gamma$ is the linear form \[ \tr_{\Nh\Gamma} : \Nh \Gamma \longrightarrow \C \] mapping $A$ to $\tr_{\Nh\Gamma} (A) = \langle A (e (1)) , e(1)\rangle $. The trace is faithful in the sense that $\tr_{\Nh\Gamma} A = 0$ for a positive operator $A$ in $\Nh\Gamma$ implies that $A = 0$. Moreover $\tr_{\Nh\Gamma}$ vanishes on commutators and satisfies the estimate $|\tr_{\Nh\Gamma} A| \le \|A\|$. On $L^1 (\Gamma ,\mathbb{C})$ it is given by $\tr_{\Nh\Gamma} (w) = w (1)$. All this is easy to check.

The Fuglede--Kadison determinant of $A$ in $(\Nh\Gamma)^*$ is the positive real number
	\begin{equation}
	\label{eq:29}
\ddet_{\Nh\Gamma} A = \exp \left( \halb \tr_{\Nh\Gamma} (\log AA^*) \right) .
	\end{equation}

Note here that the operator $AA^*$ is invertible and positive so that $\log AA^*$ is defined by the functional calculus. If $E_{\lambda}$ is a spectral resolution for $AA^*$ we can also write:
	\begin{equation}
	\label{eq:30}
\ddet_{\Nh\Gamma} A = \exp \left( \halb \int^{\infty}_0 \log \lambda \; d \, \tr_{\Nh\Gamma} (E_{\lambda}) \right) .
	\end{equation}
If the group $\Gamma$ is finite we have $\Nh\Gamma = \C \Gamma$ and
	\begin{equation}
	\label{eq:31}
\ddet_{\Nh\Gamma} A = |\ddet A|^{1 / |\Gamma|} .
	\end{equation}
We refer to \cite{L} Example 3.13 or \cite{De} Proposition 3.1 for a discussion of the case $\Gamma = \Z^n$, or see below.

For positive $A$ in $(\Nh\Gamma)^*$ we have \[ \ddet_{\Nh\Gamma} A = \exp \tr_{\Nh\Gamma} (\log A) . \] For operators $0 \le A \le B$ in $(\Nh\Gamma)^*$ this implies the inequality
	\begin{equation}
	\label{eq:32}
\ddet_{\Nh\Gamma} A \le \ddet_{\Nh\Gamma} B .
	\end{equation}
It is a non-trivial fact that $\ddet_{\Nh\Gamma}$ defines a {\it homomorphism} \[ \ddet_{\Nh\Gamma} : (\Nh\Gamma)^* \longrightarrow \R^*_+ . \] The intuitive reason for this is the Campbell--Hausdorff formula. The actual proof in \cite{FK} for $I\!I_1$-factors and the analogous argument for group von~Neumann algebras in \cite[3.2]{L} is different though in order to avoid convergence problems.

The following remark on the disintegration of $\ddet_{\Nh\Gamma}$ serves only to give some further background on this determinant. The required theory can be found in \cite{Di} for example.

Let $Z (\Gamma)$ be the center of $\Gamma$. For a character $\chi$ in $\widehat{Z (\Gamma)}$ consider the Hilbert space \[ L^2 (\Gamma ,\mathbb{C})_{\chi} = \{ h \; \mbox{in} \; L^2 (\Gamma ,\mathbb{C}) : \rho^{\gamma} h = \chi (\gamma) h \; \mbox{for all} \; \gamma \; \mbox{in} \; Z (\Gamma) \} . \] By the spectral theorem applied to the commuting unitary operators $\rho^{\gamma}$ with $\gamma$ in $Z (\Gamma)$ we have a direct integral decomposition \[ L^2 (\Gamma,\mathbb{C}) = \int_{\widehat{Z (\Gamma)}} L^2 (\Gamma,\mathbb{C})_{\chi} \, d\mu (\chi) \] where $\mu$ is the Haar probability measure on the compact abelian group $\widehat{Z (\Gamma)}$. Correspondingly we get a disintegration \[ \Nh\Gamma = \int_{\widehat{Z (\Gamma)}} \Nh_{\chi} \Gamma \, d\mu (\chi) \] where $\Nh_{\chi} \Gamma$ is the von~Neumann algebra of bounded operators on $L^2 (\Gamma,\mathbb{C})_{\chi}$ which commute with the left $\Gamma$-action $\lambda$.

The Fuglede--Kadison determinant $\ddet_{\Nh_{\chi} \Gamma}$ is defined similarly as above and we get the formula \[ \log \ddet_{\Nh\Gamma} A = \int_{\widehat{Z (\Gamma)}} \log \ddet_{\Nh_{\chi} \Gamma} A_{\chi} \, d\mu (\chi) . \] Here $A_{\chi}$ is the restriction of $A$ to an invertible operator from $L^2 (\Gamma,\mathbb{C})_{\chi}$ to itself. If $\Gamma$ is commutative we find that \[ \log \ddet_{\Nh\Gamma} A = \int_{\hat{\Gamma}} \log |a (\chi)| \, d\mu (\chi) , \] if $A_{\chi}$ acts on the $1$-dimensional space $L^2 (\Gamma,\mathbb{C})_{\chi}$ by multiplication with $a (\chi) \in \C^*$. For $f$ in $\C\Gamma$ in particular we get
	\begin{equation}
	\label{eq:33}
\log \ddet_{\Nh\Gamma} f = \int_{\hat{\Gamma}} \log |\hat{f}(\chi)| \, d \mu (\chi)
	\end{equation}
where $\hat{f}$ is the Fourier transform of $f$.

We now begin with some preparations for the proof of the entropy formula in Theorem \ref{theorem1}.

Let $\Gamma$ be a countable residually finite discrete group and let $(\Gamma_n , n \ge 1)$ be a sequence of finite index normal subgroups with $\lim_{n\to \infty} \Gamma_n = \{ 1 \}$. Set $\Gamma^{(n)} = \Gamma_n /\Gamma$. For $f$ in $L^1 (\Gamma,\mathbb{C})$ the bounded operator $\rho_f : L^2 (\Gamma,\mathbb{C}) \to L^2 (\Gamma,\mathbb{C})$ given by right convolution with $f^*$ satisfies the norm estimate
	\begin{equation}
	\label{eq:34}
\|\rho_f \| \le \| f \|_1 .
	\end{equation}
The group $\Gamma$ acts via $\lambda$ on $L^{\infty} (\Gamma,\mathbb{C})$ and we have an isomorphism of finite dimensional $\C$-vector spaces
	\begin{equation}
	\label{eq:35}
L^{\infty} (\Gamma,\mathbb{C})^{\Gamma_n} \cong L^{\infty} (\Gamma^{(n)},\mathbb{C})
	\end{equation}
given by viewing left $\Gamma_n$-invariant functions on $\Gamma$ as functions on $\Gamma_n / \Gamma$. Since $\rho_f$ is left $\Gamma_n$-equivariant it induces an endomorphism of $L^{\infty} (\Gamma,\mathbb{C})^{\Gamma_n}$ and hence an endomorphism of $L^{\infty} (\Gamma^{(n)},\mathbb{C}) = \C \Gamma^{(n)} = L^2 (\Gamma^{(n)},\mathbb{C})$ which we denote by:
	\begin{equation}
	\label{eq:36}
\rho^{(n)}_f : L^2 (\Gamma^{(n)},\mathbb{C}) \longrightarrow L^2 (\Gamma^{(n)},\mathbb{C}) .
	\end{equation}
Consider the map ``integration along the fibres'':
	\begin{equation}
	\label{eq:37}
L^1 (\Gamma,\mathbb{C}) \longrightarrow L^1 (\Gamma^{(n)},\mathbb{C})
	\end{equation}
given by sending $f : \Gamma \to \C$ to the function $f^{(n)} : \Gamma^{(n)} \to \C$ defined by \[ f^{(n)} (\delta) = \sum_{\gamma \in \delta} f (\gamma) , \] for all residue classes in $\delta$ in $\Gamma^{(n)}$. As short calculation shows that \eqref{eq:37} is a homomorphism of $\C$-algebras with involution such that $\|f^{(n)} \|_1 \le \| f \|_1$. Moreover we have
	\begin{equation}
	\label{eq:38}
\rho^{(n)}_f = \rho_{f^{(n)}} : L^2 (\Gamma^{(n)},\mathbb{C}) \longrightarrow L^2 (\Gamma^{(n)},\mathbb{C}) .
	\end{equation}
By the estimate \eqref{eq:34} applied to $f^{(n)}$ and $\Gamma^{(n)}$ we have $\|\rho_{f^{(n)}}\| \le \| f^{(n)} \|_1$. Using \eqref{eq:38} we get the uniform estimate for all $n \ge 1$
	\begin{equation}
	\label{eq:39}
\| \rho^{(n)}_f \| \le \| f \|_1 .
	\end{equation}
From the definition of $\rho^{(n)}_f$ the relation $\rho^{(n)}_{fg} = \rho^{(n)}_f \rho^{(n)}_g$ follows. Hence $\rho^{(n)}_f$ is invertible if $f$ is invertible in $L^1 (\Gamma,\mathbb{C})$ and we have $(\rho^{(n)}_f)^{-1} = \rho^{(n)}_{f^{-1}}$. Together with \eqref{eq:39} this gives the estimate
	\begin{equation}
	\label{eq:40}
\| (\rho^{(n)}_f)^{-1} \| \le \| f^{-1} \|_1 \quad \mbox{for} \; f \; \mbox{in} \; L^1 (\Gamma,\mathbb{C})^* \; \mbox{and all} \; n \ge 1 .
	\end{equation}
Note that by equation \eqref{eq:31} we have
	\begin{eqnarray}
	\label{eq:41}
\log \ddet_{\Nh\Gamma^{(n)}} f^{(n)} & = & \frac{1}{|\Gamma / \Gamma_n|} \log |\!\det \rho_{f^{(n)}}|
	\\
& = & \frac{1}{|\Gamma / \Gamma_n|} \log |\ddet (\rho_f \, |_{L^{\infty} (\Gamma,\mathbb{C})^{ \Gamma_n}})| . \nonumber
	\end{eqnarray}
Hence Theorem \ref{theorem1} of the introduction follows from Theorem \ref{t:expansive}, Theorem \ref{t:entropy} and the following result:

	\begin{theo}
	\label{theorem61}
Let $\Gamma$ be a countable residually finite discrete group and $(\Gamma_n , n \ge 1)$ a sequence of finite index normal subgroups with $\lim_{n\to \infty} \Gamma_n = \{ 1 \}$. For $f$ in $L^1 (\Gamma,\mathbb{C})^*$ we have \[ \ddet_{\Nh\Gamma} f = \lim_{n \to \infty} \ddet_{\Nh\Gamma^{(n)}} f^{(n)} . \]
	\end{theo}

	\begin{proof}
Because of the relation $(ff^*)^{(n)} = f^{(n)} f^{(n)*}$ the assertion means:
	\begin{equation}
	\label{eq:42}
\tr_{\Nh\Gamma} \log \rho_g = \lim_{n\to \infty} \tr_{\Nh\Gamma^{(n)}} \log \rho_{g^{(n)}}
	\end{equation}
for $g = ff^*$ in $L^1 (\Gamma,\mathbb{C})^*$. But $\rho_g = \rho_f \rho^*_f$ and $\rho_{g^{(n)}} = \rho_{f^{(n)}} \rho^*_{f^{(n)}}$ are positive operators on $L^2 (\Gamma,\mathbb{C})$. By the estimate \eqref{eq:34} applied to $g$ and $g^{-1}$ and the estimates \eqref{eq:39} and \eqref{eq:40} applied to $g$ instead of $f$ it follows that the spectra $\sigma (\rho_g)$ and $\sigma (\rho_{g^{(n)}})$ lie in the closed interval $I = [ \|g\|^{-1}_1 , \|g \|_1]$.

Fix $\varepsilon > 0$. By the Weierstrass approximation theorem there exists a real polynomial $Q$ such that \[ \sup_{t \in I} |\!\log t - Q (t)| \le \varepsilon . \] Since the spectra of $\rho_g$ and $\rho_{g^{(n)}}$ lie in $I$ it follows that we have \[ \| \log \rho_g - Q (\rho_g) \| \le \varepsilon \quad \mbox{and} \quad \| \log \rho_{g^{(n)}} - Q (\rho_{g^{(n)}}) \| \le \varepsilon . \] Using the estimate $|\tr_{\Nh\Gamma} A | \le \| A \|$ we find:
	\begin{eqnarray*}
\lefteqn{|\tr_{\Nh \Gamma} \log \rho_g - \tr_{\Nh\Gamma^{(n)}} \log \rho_{g^{(n)}}| \le |\tr_{\Nh\Gamma} (\log \rho_g - Q (\rho_g))| +}
	\\
& & |\tr_{\Nh\Gamma} Q (\rho_g) - \tr_{\Nh\Gamma^{(n)}} Q (\rho_{g^{(n)}})| + |\tr_{\Nh\Gamma^{(n)}} (\log \rho_{g^{(n)}} - Q (\rho_{g^{(n)}})) |
	\\
& \le & \|\log \rho_g - Q (\rho_g) \| + |\tr_{\Nh\Gamma} Q (\rho_g) - \tr_{\Nh\Gamma^{(n)}} Q (\rho_{g^{(n)}})| + \|\log \rho_{g^{(n)}} - Q (\rho_{g^{(n)}}) \|
	\\
& \le & 2 \varepsilon + |\tr_{\Nh\Gamma} Q (g) - \tr_{\Nh\Gamma^{(n)}} Q (g^{(n)}) |.
	\end{eqnarray*}
Hence formula \eqref{eq:42} is a consequence of the following Lemma \ref{lemma62}.
	\end{proof}

	\begin{lemm}
	\label{lemma62}
For any $f$ in $L^1 (\Gamma,\mathbb{C})$ and any complex polynomial $Q (t)$ we have \[ \tr_{\Nh\Gamma} Q (f) = \lim_{n\to \infty} \tr_{\Nh\Gamma^{(n)}} Q (f^{(n)}) . \]
	\end{lemm}

	\begin{proof}
Since $Q (f)$ is in $L^1 (\Gamma,\mathbb{C})$ and $Q (f)^{(n)} = Q (f^{(n)})$ it suffices to prove the assertion for $Q (t) = t$ i.e. that \[ \tr_{\Nh\Gamma} f = \lim_{n\to\infty} \tr_{\Nh\Gamma^{(n)}} f^{(n)} \] for all $f$ in $L^1 (\Gamma ,\mathbb{C})$. Writing $f = \sum_{\gamma} f_\gamma e (\gamma)$, we have $f^{(n)} = \sum_{\gamma} f_\gamma e (\overline{\gamma})$ where $\overline{\gamma} = \Gamma_n \gamma$. Thus we get \[ \tr_{\Nh\Gamma} f = f_1 \quad \mbox{and} \quad \tr_{\Nh\Gamma^{(n)}} f^{(n)} = \sum_{\gamma \in\Gamma_n} f_\gamma . \] Fix some $\varepsilon > 0$. Since $f$ is in $L^1 (\Gamma,\mathbb{C})$ we have $\sum_{\gamma \in \Gamma} |f_\gamma | < \infty$. Hence there is a finite subset $K$ of $\Gamma$ with $1 \in K$ such that we have $\sum_{\gamma \in \Gamma \smallsetminus K} |f_\gamma | < \varepsilon$. Since $\lim_{n\to \infty} \Gamma_n = \{ 1 \}$ we can find an index $N \ge 1$ such that $\Gamma_n \cap K^{-1} K = \{ 1 \}$ for all $n \ge N$. Since $1 \in K$ it follows that $\Gamma_n \cap K = \{ 1 \}$ for all $n \ge N$ as well. For $n \ge N$ we therefore get the estimate:
	\begin{eqnarray*}
|\tr_{\Nh\Gamma} f - \tr_{\Nh\Gamma^{(n)}} f^{(n)}| & = & |f_1 - \sum_{\gamma \in \Gamma_n} f_\gamma | \le \sum_{\gamma \in \Gamma_n \smallsetminus \{ 1 \} } |f_\gamma |
	\\
& \le & \sum_{\gamma \in \Gamma \smallsetminus K} |f_\gamma | < \varepsilon .
	\end{eqnarray*}
Since $\varepsilon > 0$ was arbitrary the lemma is proved.
	\end{proof}

	\begin{rema}
For $f$ in $\Z\Gamma$ the element $f^{(n)}$ lies in $\Z \Gamma^{(n)}$. Hence $\tr_{\Nh\Gamma} f$ and $\tr_{\Nh\Gamma^{(n)}} f^{(n)}$ are integers and it follows that we have $\tr_{\Nh\Gamma} f = \tr_{\Nh\Gamma^{(n)}} f^{(n)}$ for all $n \gg 0$, c.f. \cite{C} Lemma 2.2 or \cite{Schick} 5.5. Lemma.
	\end{rema}

For an element $f$ in $L^1 (\Gamma,\mathbb{C})$ let $\sigma_f = \lim_{n\to \infty} \|f^n\|^{1/n}_1 \le \|f\|_1$ be the spectral radius of $f$. We have
	\begin{eqnarray*}
\log \ddet_{\Nh\Gamma} f & = & \halb \tr_{\Nh\Gamma} \log \rho_{ff^*} \le \halb \|\log \rho_{ff^*} \| \overset{\eqref{eq:34}}{\le} \halb \log \|ff^* \|_1
	\\
& \le & \log \|f \|_1 .
	\end{eqnarray*}
Hence
	\begin{equation}
	\label{eq:43}
\log \ddet_{\Nh\Gamma} f = \frac{1}{n} \log \ddet_{\Nh\Gamma} f^n \le \log \| f^n\|^{1/n}_1 .
	\end{equation}
Therefore \[ \log \ddet_{\Nh\Gamma} f \le \log \sigma_f . \] If $f$ is in $L^1 (\Gamma,\mathbb{C})^*$ inequality \eqref{eq:43} applied to $f^{-1}$ gives: \[ \log \ddet_{\Nh\Gamma} f \ge -\log \|f^{-n}\|^{1/n}_1 \] and in the limit \[ \log \ddet_{\Nh\Gamma} f \ge - \log \sigma_{f^{-1}} . \] Together with Theorem \ref{theorem1} this gives the following corollary:

	\begin{coro}
	\label{cor63}
In the situation of Theorem \ref{theorem1}, if $f \in \Z \Gamma$ is a unit in $L^1 (\Gamma,\mathbb{C})$ we have the estimates \[ \log \|f^{-1}\|^{-1}_1 \le \log (\sigma_{f^{-1}})^{-1} \le h_f \le \log \sigma_f \le \log \|f \|_1 . \]
	\end{coro}

	\begin{exam}
For $\Gamma = \Z^n$ the spectrum of the commutative normed $*$-algebra $L^1 (\Z^n)$ is the real $n$-torus $\mathbb{T}^n\cong\mathbb{S}^n$, where $\mathbb{S}=\{z\in\mathbb{C}:|z|=1\}$ (c.f. \cite[XI 2. Example 2]{Y}). Hence an element $f$ in $L^1 (\Z^n)$ is invertible if and only if the Fourier series $\hat{f}$ does not vanish in any point of $\mathbb{S}^n $ (Wiener's theorem). By \cite{Y} XI 2. Theorem 1, we have $\sigma_f = \max_{x \in \mathbb{S}^n} |\hat{f} (x)|$ and $\sigma_{f^{-1}} = \max_{x \in \mathbb{S}^n} |\widehat{f^{-1}} (x)| = \max_{x \in \mathbb{S}^n} |\hat{f} (x)^{-1}|$. For an element $f$ in $\Z [\Z^n]$ which is invertible in $L^1 (\Z^n)$ the corollary therefore gives the estimates: \[ \min_{x \in \mathbb{S}^n} \log |\hat{f}(x)| = \log (\sigma_{f^{-1}})^{-1} \le h_f \le \log \sigma_f = \max_{x \in \mathbb{S}^n} \log |\hat{f} (x)| . \] Of course, this estimate follows directly from the formula \cite{LSW} \[ h_f = \int_{\mathbb{S}^n} \log |\hat{f} (x) | d\mu (x) , \] which in our case is also a special case of Theorem \ref{theorem1} and equation \eqref{eq:33}.
	\end{exam}

	\begin{coro}
	\label{cor64}
In the situation of Theorem \ref{theorem1}, consider elements $f,g$ in $\Z\Gamma \cap L^1 (\Gamma,\mathbb{C})^*$ which satisfy $0 \le f \le g$ in $\Nh\Gamma$. Then we have $h_f \le h_g$ and equality $h_f = h_g$ holds if and only if $f = g$.
	\end{coro}

	\begin{proof}
The first assertion follows from inequality \eqref{eq:32} applied to the operators $\rho_f$ and $\rho_g$ on $L^2 (\Gamma,\mathbb{C})$. By the monotonicity of the logarithm the operator $A = \log \rho_g - \log \rho_f$ in $\Nh\Gamma$ is positive. If the entropies are equal, i.e. if $h_f = h_g$, then it follows from Theorem \ref{theorem1} that \[ \tr_{\Nh\Gamma} (\log \rho_f) = \tr_{\Nh\Gamma} (\log \rho_g) . \] But this means that $\tr_{\Nh\Gamma} A = 0$ and hence $A = 0$ since $\tr_{\Nh\Gamma}$ is faithful. Applying $\exp$ to the equality $\log \rho_f = \log \rho_g$ gives the second assertion.
	\end{proof}

The following is an application of our theory to the structure of $(\Nh\Gamma)^*$.

	\begin{coro}
	\label{cor65}
For a group $\Gamma$ as in Theorem \ref{theorem1} let $f$ be an element of $\Z\Gamma$ which is invertible in $L^1 (\Gamma)$ but does not have a left inverse in $\Z\Gamma$. Then $f$ is not contained in the commutator subgroup of $(\Nh\Gamma)^*$.
	\end{coro}

	\begin{proof}
Otherwise we would have $\ddet_{\Nh\Gamma} f = 1$ and hence $h_f = 0$ by Theorem \ref{theorem1}. By assumption we have $\Z\Gamma \neq \Z\Gamma f$ and hence $X_f \neq 0$. But then Theorem \ref{t:positive} tells us that $h_f > 0$, contradiction.
	\end{proof}


\begin{thebibliography}{99}


\bibitem{C}B. Clair, {\it Residual amenability and the approximation of $L^2$-invariants}, Michigan Math. J. {\bf 46} (1999), 331--346. 

\bibitem{Choi} Y. Choi, \textit{Injective convolution operators on $L^\infty(\Gamma)$ are surjective}, arXiv:math.FA/0606367 v1 15 Jun 2006.

\bibitem{De}C. Deninger, {\it Fuglede--Kadison determinants and entropy for actions of discrete amenable groups}, J. Amer. Math. Soc. \textbf{19} (2006), 737--758. 

\bibitem{Di}J. Dixmier, {\it Von~Neumann algebras}, North Holland, Amsterdam, 1981. 

\bibitem{DS} N. Dunford and J.T. Schwartz, \textit{Linear operators}, vol. 1, Wiley, New York, 1967. 

\bibitem{ER} M. Einsiedler and H. Rindler, \textit{Algebraic actions of the discrete Heisenberg group and other non-abelian groups}, Aequationes Math. \textbf{62} (2001), 117--135.

\bibitem{FK}B. Fuglede and R.V. Kadison, {\it Determinant theory in finite factors}, Ann. Math. {\bf 55} (1952), 520--530. 

\bibitem{Kap} I. Kaplansky, \textit{Fields and rings}, University of Chicago Press, Chicago-London, 1969.

\bibitem{LS} D. Lind and K. Schmidt, \textit{Homoclinic points of algebraic $\mathbf{Z}^d$-actions}, J. Amer. Math. Soc. \textbf{12} (1999), 953--980. 

\bibitem{LSW} D. Lind, K. Schmidt and T. Ward, {\it Mahler measure and entropy for commuting automorphisms of compact groups}, Invent. Math. {\bf 101} (1990), 593--629. 

\bibitem{L} W. L\"uck, {\it $L^2$-invariants: Theory and applications to geometry and $K$-theory}, Springer, Berlin, 2002. 



\bibitem{Schick}T. Schick, {\it $L^2$-determinant class and approximation of $L^2$-Betti numbers}, TAMS {\bf 353} (2001), 3247--3265.

\bibitem{Sch}K. Schmidt, {\it Dynamical systems of algebraic origin}, Birkh\"auser, Basel 1995. 

\bibitem{Weiss} B. Weiss, \textit{Monotileable amenable groups}, Amer. Math. Soc. Transl. \textbf{202} (2001), 257--262.


\bibitem{Y} K. Yosida, {\it Functional analysis}, Springer, Heidelberg, 1971.
	\end{thebibliography}
	\end{document}